\documentclass[a4paper,12pt]{article}
\usepackage{amsmath}
\usepackage{amssymb}
\usepackage{tabularx}
\usepackage{enumerate}
\usepackage{graphicx}

\newtheorem{thm}{Theorem}[section]
\newtheorem{lem}[thm]{Lemma}
\newtheorem{cor}[thm]{Corollary}

\newtheorem{rmk}{Remark}[section]
\newtheorem{defi}{Definition}[section]
\newtheorem{pppp}{Proof}

\newcommand{\qed}{\hspace{1em}\mbox{\raisebox{0.65ex}{\fbox{}}}}

\numberwithin{equation}{section}

\newcommand{\be}{\begin{equation}}
\newcommand{\ee}{\end{equation}}
\newcommand\bes{\begin{eqnarray}} \newcommand\ees{\end{eqnarray}}
\newcommand{\bess}{\begin{eqnarray*}}
\newcommand{\eess}{\end{eqnarray*}}

\newcommand{\R}{\mathbb{R}}
\newcommand{\bpf}{{\bf Proof:\ \ }}
\newcommand{\epf}{\mbox{}\hfill $\Box$}

\begin{document}

\thispagestyle{empty}

\title{The invasive dynamics of Aedes aegypti mosquito in a heterogenous environment\thanks{The work is partially supported by the NSFC of China (Grant No.11371311, 11501494, 11571301), the High-End Talent Plan of Yangzhou University, China.}}
\date{\empty}

\author{Mengyun Zhang$^{a}$, Jing Ge$^{a, b}$ and Zhigui Lin$^{a}\thanks{Corresponding author. Email: zglin68@hotmail.com (Z. Lin).}$\\
{\small $^a$ School of Mathematical Science, Yangzhou University, Yangzhou 225002, China}\\
{\small $^b$ School of Mathematical Science, Huaiyin Normal University, Huaian 223300, China}}
 \maketitle

\begin{quote}
\noindent
{\bf Abstract.} { 
\small A reaction-diffusion-advection model is proposed and investigated to understand the invasive dynamics of Aedes aegypti mosquitoes. The free boundary is introduced to
model the expanding front of the invasive mosquitoes in a heterogenous environment.
The threshold $R^D_0$ for the model with Dirichlet boundary condition is defined and the threshold $R^F_0(t)$ for the free boundary problem is introduced, and the long-time behavior of positive solutions to the reaction-diffusion-advection system is discussed.
Sufficient conditions for the mosquitoes to be eradicated or to spread are given. We show that, if $R^F_0(\infty)\leq 1$,
 the mosquitoes always vanish, and if $R^F_0(t_0)\geq 1$ for some $t_0\geq 0$, the mosquitoes
must spread, while if $R^F_0(0)<1<R^F_0(\infty)$, the spreading or vanishing of the mosquitoes depends on the initial number of  mosquitoes,
or mosquitoes' invasive ability on the free boundary. Moreover, numerical simulations indicate that the advection and the expanding capability affect the mosquitoes' invasive fronts . }

\noindent {\it MSC:} primary: 35R35; secondary: 35K60

\medskip
\noindent {\it Keywords: } Invasive mosquitoes; Reaction-diffusion-advection model; Free boundary; Spreading and vanishing

\end{quote}

\section{Introduction}
Mosquitoes cause more human suffering than any other organism -- over one million people worldwide die from mosquito-borne diseases every year. The Aedes mosquitoes are responsible for transmitting several of the most debilitating mosquito-borne viruses, among which Aedes aegypti is the primary transmitter. These include dengue fever \cite{WHO}, Zika \cite{AFDM,DGYL}, chikungunya \cite{SMC}, etc.

There are nearly 400 million people infected by dengue fever each year, and as a result an estimated 25,000 deaths. And for Zika virus (ZIKV)  \cite{DGYL}, which is  closely related to dengue, and primarily transmitted to humans through the bites of infected female mosquitoes from the Aedes genus, almost one in five infected individuals develops symptoms like  rash, conjunctivitis,  mild fever and joint pain. Since no effective treatment or vaccine is available to treat or prevent these mosquito-borne diseases currently, and an infected mosquito is able to transmit virus in its remaining life, an effective method for the mosquito-borne disease control is vector eradication.

 Aedes aegypti bites primarily during the day, both indoors and outdoors, which acts most
actively for approximately two hours after sunrise and several hours before sunset, but can
also bite at night in well-lit areas. This mosquito can bite people without being noticed, and prefers biting people but also bites dogs and other domestic animals, mainly mammals. Thus the mosquito could easily transmit viruses mentioned above.
Aedes aegypti control aims to reduce the density of adult mosquitoes populations under a threshold where the Aedes aegypti-borne epidemics could not occur.

As we know, the evolution process of the habitat plays an important role in investigating the dynamics of invasive species, and to describe this, the free boundary problems have been studied in many areas \cite{CS}. The well-known Stefan condition has been applied to describe the interaction and expanding process at the boundary, for example, the spreading of the invasive populations in \cite{DPW, DGP, DL, GW, KY, Wa, WZ, ZX}.

To explore the temporal and spatial dispersal of the Aedes aegypti mosquitoes, Tian and Ruan \cite{CTSR} proposed an advection-reaction-diffusion model with free boundary, based on \cite{LNW}, where the vector mosquitoes population is divided into two life stages: the winged form (adult female mosquitoes) and an aquatic population (including eggs, larvae and pupae):
\begin{eqnarray}
\left\{
\begin{array}{ll}
M_t=DM_{xx}-\nu M_x+\gamma A(1-\frac{M}{K_1})-\mu_1M,&t>0, 0<x<h(t), \\
M(t,x)=0, &t>0, x\geq h(t), \\
A_t=r(1-\frac{A}{K_2})M-(\mu_2+\gamma)A,&t>0, x>0, \\
M_x(t,0)=0, A_x(t,0)=0, &t>0,\\
M(t,h(t))=0, h'(t)=-\mu M_x(t,h(t)),&t>0,\\
h(0)=h_0,M(0,x)=M_0(x), &x \in [0,h_0],\\
A(0,x)=A_0(x),&x\in [0,\infty),
\end{array} \right.
\label{Aa01}
\end{eqnarray}
the free boundary $h(t)$ describes the expanding front of mosquitoes, which focuses on the changing of invasive habitat, $M(t,x)$ and $A(t,x)$ represent the density of the winged and aquatic mosquitoes at time $t$ and location $x$, respectively, $D$ denotes the diffusion rate as a result of the random walk of the winged mosquitoes, $\nu$ is the advection rate caused by wind, $r$ is the oviposition rate of adult female mosquitoes, while $\gamma$ is the rate of maturation from aquatic form into the winged mosquitoes. $K_1$ and $K_2$ are the carrying capacity of winged and aquatic mosquitoes, respectively. $\mu_1$ and $\mu_2$ denote, respectively, the death rates of  the winged mosquitoes and the aquatic form.
They gave sufficient conditions for mosquitoes to be spreading or vanishing.

As we can see in (\ref{Aa01}), the environment is assumed to be homogenous, 
however, environmental heterogeneity has been recognized as more and more important factor to the persistence of infectious diseases or the spreading of the invasive species \cite{AL, JKZH, LPW, ZX}, in particular the authors in \cite{ZX} divided the environment into two cases: strong heterogeneous environment and weak heterogeneous environment. Recently, Allen, et al \cite{AL} proposed an SIS reaction-diffusion model to study the dynamics of the transmission of infectious diseases in a heterogeneous environment:
\begin{eqnarray}
\left\{
\begin{array}{lll}
S_{t}-d_S\Delta S=-\frac{\beta (x) SI}{S+I}+\gamma(x)I,\; &  x\in\Omega, \ t>0, \\
I_{t}-d_I\Delta I=\frac{\beta (x) SI}{S+I}-\gamma(x)I,\; &   x\in \Omega,\ t>0,
\end{array} \right.
\label{Aa02}
\end{eqnarray}
where $S(x,t)$ and $I(x,t)$ denote respectively, the density of susceptible and infected individuals at location $x$ and time $t$. $d_S$ and $d_I$ denote the positive diffusion rates for the susceptible and infected individuals. $\beta (x)$ and $\gamma(x)$ account for the contact transmission rate and recovery rate of the disease at $x$, which are spatial-dependent, respectively. Their results show that spatial heterogeneity has great influence on the persistence and extinction of the disease.

Considering the spatial heterogeneity, based on model (\ref{Aa01}), we consider a reaction-diffusion-advection problem with free boundaries $x=g(t)$ and $x=h(t)$ to describe the spatial dispersal dynamics of A. aegypti mosquitoes:
{\small \begin{eqnarray}
\left\{
\begin{array}{ll}
M_{t}=D M_{xx}-\nu M_{x}+\gamma (x)A(1-\frac{M}{K_1}) -\mu_1(x)M,&t>0,\; g(t)<x<h(t),   \\
A_{t}=r (x)(1-\frac{A}{K_2})M -(\mu_2(x)+\gamma(x))A,\; &t>0,\; g(t)<x<h(t),   \\
M(t,x)=A(t,x)=0,\, &t>0,  x=g(t)\, \textrm{or}\, x=h(t),\\
g(0)=-h_0,\; g'(t)=-\mu M_{x}(t, g(t)), & t>0, \\
 h(0)=h_0, \; h'(t)=-\mu M_{x}(t, h(t)), & t>0,\\
M(0,x)=M_0(x),\ A(0,x)=A_0(x),&-h_0\leq x\leq h_0,
\end{array} \right.
\label{a3}
\end{eqnarray}}
where $\gamma (x)$, $r (x)$, $\mu_1(x)$ and $\mu_2(x)$ are corresponding spatial-dependent rates to model (\ref{Aa01}), $x=g(t)$ and $x=h(t)$ are the moving left and right
boundaries to be determined, $h_0$ and $\mu $ are positive constants, and the initial functions
$M_{0}$ and $A_0$ are nonnegative and satisfy
\begin{eqnarray}
\left\{
\begin{array}{lll}
M_{0}\in C^2([-h_0, h_0]),\, M_{0}(\pm h_0)=0\, \textrm{and} \, 0< M_{0}(x)<K_1,\,  x\in (-h_0, h_0), \\
A_{0}\in C^2([-h_0, h_0]), A_{0}(\pm h_0)=0\, \textrm{and} \, 0< A_{0}(x)<K_2,\,  x\in (-h_0, h_0).
\end{array} \right.
\label{Ae1}
\end{eqnarray}
Furthermore, we assume
$$(H)\;\;\;\;\lim_{x\rightarrow\pm \infty}r(x)=r_{\infty},\;\;\;
\lim_{x\rightarrow \pm \infty}\gamma(x)=\gamma_{\infty},\;\;\;
\lim_{x\rightarrow \pm \infty}\mu_1(x)={\mu_1}_{\infty},\;\;\;$$
$$\lim_{x\rightarrow \pm \infty}\mu_2(x)={\mu_2}_{\infty}
\ \textrm{and}\ \frac{r_{\infty}\gamma_{\infty}}{{\mu_2}_{\infty}+\gamma_{\infty}}-{\mu_1}_{\infty}>0,$$
 which implies that far sites of the habitat
are similar and high-risk.
And we assume that $\nu<2D\sqrt{\frac{r_\infty\gamma_\infty-{\mu_1}_\infty({\mu_2}_\infty+\gamma_\infty)}{{\mu_2}_\infty+\gamma_\infty}}$,
which represents the small advection.

This paper is organized as follows.
In Section 2, we present the global existence and uniqueness of the solution to
(\ref{a3}) by applying the contraction mapping theorem, comparison principle is also presented. Section 3 deals with some thresholds and the related properties. Section 4 is devoted to sufficient conditions for mosquitoes to vanish. The case and conditions for mosquitoes to spread are discussed in Section 5. The paper ends with 
 a brief discussion.

\section{Existence and uniqueness}
First, we present the following local existence and
uniqueness results by using the contraction mapping theorem, and then we show global existence with the help of using
some suitable estimates.
\begin{thm} For any given $(M_0, A_0)$ satisfying \eqref{Ae1},  and any $\alpha \in (0, 1)$, there exists a $T>0$ such that
problem  \eqref{a3} admits a unique solution
$$(M, A; g, h)\in [C^{(1+\alpha)/2, 1+\alpha}(D_{T})]^2\times [C^{1+\alpha/2}([0,T])]^2;$$
moreover,
{\small \begin{eqnarray} \|M\|_{C^{(1+\alpha)/2,1+\alpha}({D}_{T})}+\|A\|_{C^{(1+\alpha)/2,1+\alpha
}({D}_{T})}+||g\|_{C^{1+\alpha/2}([0,T])}+||h\|_{C^{1+\alpha/2}([0,T])}\leq
C,\ \label{b12}
\end{eqnarray}}
where $D_{T}=\{(t, x)\in \R^2:  t\in [0,T],x\in [g(t), h(t)]\}$, $C$
and $T$  depend only on $h_0, \alpha, \|M_0\|_{C^{2}([-h_0, h_0])}$ and $\|A_0\|_{C^{2}([-h_0, h_0])}$.
\end{thm}
\bpf
For any given $T>0$, we define
$$\mathcal{G}_{T}=\{g\in C^1([0,T]): \,  g(0)=-h_0,\,  g'(t)\leq 0,\, 0\leq t\leq T\},$$
$$\mathcal{H}_{T}=\{h\in C^1([0,T]): \, h(0)=h_0,\,  h'(t)\geq 0,\, 0\leq t\leq T\}.$$
Note that the second equation of (\ref{a3}) for  $A$ has no diffusion term, we can use $g, h$ and $M$ to represent $A$. If $g(t)\in \mathcal{G}_{T}$, $h(t)\in \mathcal{H}_{T}$ and $M(t,x)\in C(D_T)$, then for $(t,x)\in D_T$,
the unknown $A$ can be represented as
\begin{eqnarray*}
A(t,x):=H(t,M(t,x))=e^{-\frac{r(x)}{K_2}\int_0^t M(s,x)ds-\mu_2(x)t-\gamma(x)t}A(0,x)\\
+\int_0^tr(x)M(\tau,x)e^{\frac{r(x)}{K_2}\int_t^\tau M(s,x)ds+\mu_2(x)(\tau-t)+\gamma(x)(\tau-t)}d\tau.
\end{eqnarray*}

Consider the transformation $w(t,y)=M(t,x)$, $\gamma_1(t,y)=\gamma(x)$, $\mu_{11}(t,y)=\mu_1(x),$ where
$$y=\frac{2h_0x}{h(t)-g(t)}-\frac{h_0(h(t)+g(t))}{h(t)-g(t)},$$
then problem (\ref{a3}) can be transformed into
\begin{eqnarray}
\left\{
\begin{array}{ll}
w_{t}=Aw_{y}+Bw_{yy}+\gamma_1(t,y)H(t,w)(1-\frac{w}{K_1})&\\
\qquad -\mu_{11}(t,y)w,\; &t>0, \ -h_0<y<h_0, \\
w=0,\quad h'(t)=-\frac{2h_0\mu}{h(t)-g(t)}\frac{\partial w}{\partial
y},\quad &t>0, \ y=h_0,\\
w=0,\quad g'(t)=-\frac{2h_0\mu}{h(t)-g(t)}\frac{\partial w}{\partial
y},\quad &t>0, \ y=-h_0,\\
h(0)=h_0, \quad g(0)=-h_0, &\\
w(0,y)=w_0(y):=M_{0}(y), \; &-h_0\leq y\leq
h_0,
\end{array} \right.
\label{Hb}
\end{eqnarray}
where $A=A(h, g,
y)=y\frac{h'(t)-g'(t)}{h(t)-g(t)}+h_0\frac{h'(t)+g'(t)}{h(t)-g(t)}-\nu\frac{2h_0}{h(t)-g(t)}$, and
$B=B(h, g)=\frac{4h_0^2D}{(h(t)-g(t))^2}$. This transformation
changes the free boundary problem (\ref{a3}) to the initial boundary problem (\ref{Hb}) in $(-h_0, h_0)$ with more complex equations.

Similarly as those in \cite{ABL, DL}, the rest of the proof follows from the contraction mapping theorem together with the standard $L^p$ theory and the Sobolev imbedding theorem, we omit it here.
\epf

\begin{thm} Let $(M, A; g, h)$ be a solution to problem \eqref{a3} defined for $t\in (0,T_0]$ for some $T_0\in (0, +\infty)$. Then the following conclusion hold.

$(a)$ $ 0<M(t, x)\leq K_1$ and $0<A(t, x)\leq K_2 \; \mbox{ for }\, t\in (0, T_0],\; g(t)<x<h(t); $

$(b)$ There exists a constant $C_1$  independent of $T_0$ such
that
\[
 0<-g'(t),\ h'(t)\leq C_1 \; \mbox{ for } \; t\in (0,T_0]; \]

$(c)$ the solution of  \eqref{a3} exists and is unique for all $t\in(0,\infty).$
\end{thm}
\bpf
$(a)$ is directly from the comparison principle, see Lemma 2.2 in \cite{ABL}. The proof of $(b)$ is similar to that of Lemma 2.2 in \cite{DL}, where 
$$C_1:=\max \{\frac 1{2h_0},\, \frac {\nu}D+\sqrt{\frac {K_1}{2D}},\, \frac {4||M_0||_{C^1([-h_0, h_0])}}{3K_1} \} .$$
For $(c)$, since $M, A$ and $g'(t), h'(t)$ are bounded in $(g(t),h(t))\times (0, T_0]$ by constants independent of $T_0$, the maximal existing time of the solution of  \eqref{a3} can be extended to infinity.
\epf

In what follows, we exhibit the comparison principle.
\begin{lem} (The Comparison Principle)
  Assume that $\overline g, \overline
h\in C^1([0, +\infty))$, $\overline M(t, x)$, $\overline A(t, x)$ $\in C([0, +\infty)\times [\overline g(t), \overline h(t)])\cap
C^{1,2}((0, +\infty)\times(\overline g(t), \overline h(t)))$, and
\begin{eqnarray*}
\left\{
\begin{array}{lll}
\overline M_t\geq D\overline M_{xx}-\nu\overline M_x+\gamma (x)\overline A(1-\frac{\overline M}{K_1}) -\mu_1(x)\overline M,&t>0,\,\overline g(t)<x<\overline h(t),  \\
\overline A_t\geq r(x)(1-\frac{\overline A}{K_2})\overline M-(\mu_2(x)+\gamma(x))\overline A,&t>0,\,\overline g(t)<x<\overline h(t),  \\
\overline M(t, x)=\overline A(t, x)=0,\, &  t>0,\,x=\overline g(t)\textrm{or}\, x= \overline h(t),\\
\overline g(0)\leq -h_0,\; \overline g'(t)\leq -\mu \overline M_x(t, \overline g(t)), & t>0, \\
 \overline h(0)\geq h_0, \; \overline h'(t)\geq -\mu \overline M_x(t, \overline h(t)), & t>0,\\
\overline M(0,x)\geq M_0(x),\ \overline A(0,x)\geq A_0(x),&-h_0\leq x\leq h_0.
\end{array} \right.
\label{cp}
\end{eqnarray*}
Then the solution $(M, A; g, h)$ to the free boundary problem $(\ref{a3})$ satisfies
$$h(t)\leq\overline h(t),\ g(t)\geq \overline g(t),\quad t\in [0, +\infty),$$
$$M(t, x)\leq \overline M(t, x),\ A(t, x)\leq \overline A(t, x),\ t\geq 0,\,x\in [g(t), h(t)].$$\label{Com}
\end{lem}

The pair $(\overline M, \overline A; \overline g, \overline h)$ in Lemma \ref{cp} is usually called an upper solution
 to problem \eqref{a3} and a lower solution $(\underline M, \underline A; \underline g, \underline h)$ can be defined similarly by reversing all of the 
 inequalities in the obvious places.
To emphasize the dependence of the solution on the expanding capability $\mu$,
 we rewrite the solution as $(M^{\mu}, A^{\mu}; g^{\mu}, h^{\mu})$. As a corollary of Lemma \ref{cp}, we have the following monotonicity:
\begin{cor} For fixed $M_0, A_0, \nu, h_0, r(x),\gamma(x), \mu_1(x)$ and $\mu_2(x)$.
If $\mu_1\leq \mu_2$. Then $M^{\mu_1}(x, t)\leq M^{\mu_2}(x, t)$, $A^{\mu_1}(x, t)\leq A^{\mu_2}(x, t)$ in $(0, \infty)\times [g^{\mu_1}(t), h^{\mu_1}(t)]$
 and $g^{\mu_2}(t)\leq g^{\mu_1}(t)$, $h^{\mu_1}(t)\leq h^{\mu_2}(t)$ in $(0, \infty)$.
\end{cor}
\section{The threshold value}

In this section, we will give a threshold value
for the free boundary problem \eqref{a3}, which is similar to the basic reproduction number in epidemiology. First, we define a threshold value and present its properties and implications
for the following reaction-diffusion-advection model with Dirichlet boundary condition
\begin{eqnarray}
\left\{
\begin{array}{lll}
M_{t}=D M_{xx}-\nu M_{x}+\gamma (x)A(1-\frac{M}{K_1}) -\mu_1(x)M, &t>0,\;x\in (p, q),   \\
A_{t}=r (x)(1-\frac{A}{K_2})M -(\mu_2(x)+\gamma(x))A,\; &t>0,\; x\in (p, q),   \\
M(t,x)=A(t,x)=0,\, &t> 0,\,  x=p\, \textrm{or}\, x=q.
\end{array} \right.
\label{Aa2ad}
\end{eqnarray}
We linearize \eqref{Aa2ad} around (0,0) to obtain
\begin{eqnarray}
\left\{
\begin{array}{lll}
\eta_t-D\eta_{xx}+\nu\eta_x=\gamma(x)\xi-\mu_1(x)\eta,\; &t>0,\;
x\in (p, q),  \\
\xi_t=r(x) \eta-(\mu_2(x)+\gamma(x))\xi,\; &t>0,\;
x\in (p, q),  \\
\eta(t,x)=\xi(t,x)=0,&t>0,\;x=p\, \textrm{or}\, x=q,

\end{array} \right.
\label{B10f}
\end{eqnarray}
and consider the following eigenvalue problem
\begin{eqnarray}
\left\{
\begin{array}{lll}
-D\phi_{xx}+\nu\phi_x=\frac{\gamma(x)}{R_0^{DA}} \psi-\mu_1(x)\phi,\; &
x\in (p, q),  \\
0=\frac{r(x)}{R_0^{DA}} \phi-(\mu_2(x)+\gamma(x))\psi,\; &
x\in (p, q),  \\
\phi(x)=\psi(x)=0,&x=p\, \textrm{or}\, x=q.
\end{array} \right.
\label{B101f}
\end{eqnarray}
As stated in \cite{AL,CC}, we introduce the threshold value $R_0^{DA}$ by
$$R_0^{DA}=R_0^{DA}((p, q), D, \nu)=\ \sup_{\psi\in H^1_0(p, q),\psi\neq 0} \{\sqrt{\frac{\int_p^q \frac{r(x)\gamma(x)}{(\mu_2(x)+\gamma(x))}
\psi^2dx}{\int_p^q (D\psi_x^2+\frac{\nu^2}{4D}\psi^2+\mu_1(x) \psi^2)dx}}\}.$$

The following result follows from variational methods, see Chapter 2 in \cite{CC} for example.
\begin{lem} $1-R_0^{DA}$ has the same sign as $\lambda_0$, where $\lambda_0$ is the principal eigenvalue of the reaction-diffusion-advection problem
\begin{eqnarray}
\left\{
\begin{array}{lll}
-D\phi_{xx}+\nu\phi_x=\gamma(x) \psi-\mu_1(x)\phi+\lambda_0 \phi,\; &
x\in (p, q),  \\
0=r(x) \phi-(\mu_2(x)+\gamma(x))\psi+\lambda_0 \psi,\; &
x\in (p, q),  \\
\phi(x)=\psi(x)=0,&x=p\, \textrm{or}\, x=q.

\end{array} \right.
\label{B11f}
\end{eqnarray}
\end{lem}
\bpf
Substituting the second equation in \eqref{B11f} to the first equation yields
\begin{eqnarray}
\begin{array}{lll}
-D\phi_{xx}+\nu\phi_x=\frac{r(x)\gamma(x)}{\mu_2(x)+\gamma(x)-\lambda_0}\phi-\mu_1(x)\phi+\lambda_0 \phi
\end{array}.
\label{Bd1f}
\end{eqnarray}
Let $(R_0^{DA},\phi^*,\psi^*)$ be the eigen-pair of problem (\ref{B101f}), that is
\begin{eqnarray}
\left\{
\begin{array}{lll}
-D\phi^*_{xx}+\nu\phi^*_x=\frac{\gamma(x)}{R_0^{DA}} \psi^*-\mu_1(x)\phi^*,\; &
x\in (p, q),  \\
0=\frac{r(x)}{R_0^{DA}} \phi^*-(\mu_2(x)+\gamma(x))\psi^*,\; &
x\in (p, q),  \\
\phi^*(x)=\psi^*(x)=0,&x=p\, \textrm{or}\, x=q,
\end{array} \right.
\label{B12f}
\end{eqnarray}
which reduces to
\begin{eqnarray}
\begin{array}{lll}
-D\phi^*_{xx}+\nu\phi^*_x=\frac{r(x)\gamma(x)}{(R_0^{DA})^2(\mu_2(x)+\gamma(x))}\phi^*-\mu_1(x)\phi^*
\end{array}.
\label{Bd2f}
\end{eqnarray}
For convenience, taking $\Psi=e^{-\frac{\nu}{2D}x}\phi$ in (\ref{Bd1f})and $\Psi^*=e^{-\frac{\nu}{2D}x}\phi^*$ in (\ref{Bd2f}) yields
\begin{eqnarray}
\begin{array}{lll}
-D\Psi_{xx}=\frac{r(x)\gamma(x)}{\mu_2(x)+\gamma(x)-\lambda_0}\Psi-\mu_1(x)\Psi-\frac{\nu^2}{4D}\Psi+\lambda_0 \Psi
\end{array},
\label{Bd3f}
\end{eqnarray}
and
\begin{eqnarray}
\begin{array}{lll}
-D\Psi^*_{xx}=\frac{r(x)\gamma(x)}{({R_0^{DA}})^2(\mu_2(x)+\gamma(x))}\Psi^*-\mu_1(x)\Psi^*-\frac{\nu^2}{4D}\Psi^*
\end{array}.
\label{B14f}
\end{eqnarray}
By the multiply-multiply-subtract-integrate technique, we obtain
$$\int_p^q \frac{r(x)\gamma(x)}{({R_0^{DA}})^2(\mu_2(x)+\gamma(x))}\Psi\Psi^*dx=\int_p^q [\frac{r(x)\gamma(x)}{\mu_2(x)+\gamma(x)-\lambda_0}+\lambda_0]\Psi\Psi^*dx,$$
which means that
$${\textrm sign}\,(1-R_0^{DA})={\textrm sign}\,\lambda_0.$$
\epf

With the above definition, we have the following statements.
\begin{thm} The following assertions hold.

$(a)$ $R_0^{DA}$ is a positive and monotone decreasing function of $\nu$;

$(b)$ $R_0^{DA}\to 0$ as  $D\to \infty$;

$(c)$ If $\Omega_1\subseteqq \Omega_2 \subseteqq R^1$, then $R_0^{DA}(\Omega_1)\leq R_0^{DA}(\Omega_2)$, with strict inequality if $\Omega_2 \setminus \Omega_1$ is an open set. Moreover, $\lim_{(q-p)\to \infty}\, R_0^{DA}((p, q), D,\nu)\geq \sqrt{\frac {\frac{r_{\infty}\gamma_{\infty}}{{\mu_2}_{\infty}+\gamma_{\infty}}}{\frac{\nu^2}{4D^2}+{\mu_1}_{\infty}}}$ provided that $(H)$ holds;

$(d)$ If $r (x)\equiv r^*$, $\gamma(x)\equiv \gamma^*$, $\mu_1(x)\equiv \mu_1^*$, $\mu_2(x)\equiv \mu_2^*$, then
$$R_0^{DA}=\sqrt{\frac {\frac{r^*\gamma^*}{\mu_2^*+\gamma^*}}{D{(\frac{\pi}{q-p})}^2+\frac{\nu^2}{4D^2}+\mu_1^*}}.$$
\label{basic2}
\end{thm}
\bpf
The proof of part (a) is from the definition of $R_0^{DA}$, and part (d) can be obtained through direct calculations, where${(\frac{\pi}{q-p})}^2$ is the principal eigenvalue of $-\Delta$ operator with the null Dirichlet boundary condition in $(p,q)$.

For part (b), by the definition of $R_0^{DA}$ and Poinc$\acute{a}$re's inequality, we have
\begin{eqnarray*}
R_0^{DA}&=&\ \sup_{\psi\in H^1_0(p, q),\psi\neq 0} \Big\{\sqrt{\frac{\int_p^q \frac{r(x)\gamma(x)}{(\mu_2(x)+\gamma(x))}
\psi^2dx}{\int_p^q (D\psi_x^2+\frac{\nu^2}{4D}\psi^2+\mu_1(x) \psi^2)dx}}\Big\}\\
&\leq&\ \sup_{\psi\in H^1_0(p, q),\psi\neq 0} \Big\{\sqrt{\frac{\int_p^q \frac{r^M\gamma^M}{(\mu_2^m+\gamma^m)}
\psi^2dx}{\int_p^q (D\psi_x^2+\frac{\nu^2}{4D}\psi^2+\mu_1^m \psi^2)dx}}\Big\}\\
&\leq&\ \sup_{\psi\in H^1_0(p, q),\psi\neq 0} \Big\{\sqrt{\frac{\int_p^q \frac{r^M\gamma^M}{(\mu_2^m+\gamma^m)}
\psi^2dx}{\int_p^q ((D{(\frac{\pi}{q-p})}^2+\frac{\nu^2}{4D}+\mu_1^m) \psi^2)dx}}\Big\}\\
&=&\sqrt{ \frac{\frac{r^M\gamma^M}{(\mu_2^m+\gamma^m)}
} {D{(\frac{\pi}{q-p})}^2+\frac{\nu^2}{4D}+\mu_1^m}}\rightarrow 0\ \mbox{as}\ D\rightarrow\infty,
\end{eqnarray*}
where $f^M=\sup_{x\in (p,q)} \{f(x)\}, f^m=\inf_{x\in (p,q)} \{f(x)\}$ for any bounded function $f$ in $(p,q)$.

The proof of the monotonicity in (c) is similar to that of Corollary 2.3 in \cite{CC}. For the limit part, it follows from the assumption (H) that for any $\varepsilon >0$,
there exists $L_0>0$, when $|x|>L_0$, $r_\infty-\varepsilon\leq r(x)\leq r_\infty+\varepsilon, \gamma_\infty-\varepsilon\leq \gamma(x)\leq\gamma_\infty+\varepsilon, {\mu_1}_\infty-\varepsilon\leq  {\mu_1}(x)\leq  {\mu_1}_\infty+\varepsilon, {\mu_2}_\infty-\varepsilon\leq  {\mu_2}(x)\leq  {\mu_2}_\infty+\varepsilon.$

For the case $q\geq2L_0$,
\begin{eqnarray*}
& &\lim_{(q-p)\to \infty}\,R_0^{DA}((p, q), D, \nu)\\
&\geq&\ \sup_{\psi\in H^1_0(L_0, 2L_0),\psi\neq 0} \{\sqrt{\frac{\int_{L_0}^{2L_0} \frac{r(x)\gamma(x)}{(\mu_2(x)+\gamma(x))}
\psi^2dx}{\int_{L_0}^{2L_0} (D\psi_x^2+\frac{\nu^2}{4D}\psi^2+\mu_1(x) \psi^2)dx}}\}\\
&\geq&\ \sup_{\psi\in H^1_0(L_0, 2L_0),\psi\neq 0} \{\sqrt{\frac{\int_{L_0}^{2L_0} \frac{(r_\infty-\varepsilon)(\gamma_\infty-\varepsilon)}{({\mu_2}_\infty+\varepsilon)+\gamma_\infty+\varepsilon}
\psi^2dx}{\int_{L_0}^{2L_0} (D\psi_x^2+\frac{\nu^2}{4D}\psi^2+({\mu_1}_\infty+\varepsilon) \psi^2)dx}}\}\\
&\geq&\sqrt{\frac{\frac{(r_\infty-\varepsilon)(\gamma_\infty-\varepsilon)}{({\mu_2}_\infty+\varepsilon)+\gamma_\infty+\varepsilon}}
{(D{(\frac{\pi}{L_0})}^2+\frac{\nu^2}{4D}+({\mu_1}_\infty+\varepsilon)}}.
\end{eqnarray*}
Since $\varepsilon$ is arbitrary, letting $L_0\rightarrow\infty$ yields
$$\lim_{(q-p)\to \infty}\,R_0^{DA}((p, q), D, \nu)\geq \sqrt{\frac{\frac{r_\infty\gamma_\infty}{{\mu_2}_\infty+\gamma_\infty}}{{\mu_1}_\infty+\frac{\nu^2}{4D^2}}}.$$
Similarly, for $p\leq-2L_0$, we obtain the same result by replacing $(L_0,2L_0)$ with $(-2L_0,-L_0)$.
\epf

\bigskip
Noticing that the domain $(g(t), h(t))$ is changing with $t$, so the threshold value for the free boundary problem \eqref{a3} will not be a constant and
should be changing with $t$. As a result, we introduce the threshold value $R_0^F(t)$ by
\begin{eqnarray*}
R_0^F(t)&:=&R_0^{DA}((g(t),h(t))), D, \nu)\\
&=&\ \sup_{\psi\in H^1_0(g(t), h(t)),\psi\neq 0} \{\sqrt{\frac{\int_{g(t)}^{h(t)} \frac{r(x)\gamma(x)}{(\mu_2(x)+\gamma(x))}
\psi^2dx}{\int_{g(t)}^{h(t)} (D\psi_x^2+\frac{\nu^2}{4D}\psi^2+\mu_1(x) \psi^2)dx}}\}.
\end{eqnarray*}

Lemma 3.1 together with the above definition shows that
\begin{lem} $1-R_0^F(t)$ has the same sign as $\lambda_0$, where $\lambda_0$ is the principal eigenvalue of the problem
\begin{eqnarray}
\left\{
\begin{array}{lll}
-D\phi_{xx}+\nu\phi_x=\gamma(x) \psi-\mu_1(x)\phi+\lambda \phi,\; &
x\in (g(t), h(t)),  \\
0=r(x) \phi-(\mu_2(x)+\gamma(x))\psi+\lambda \psi,\; &
x\in (g(t), h(t)),  \\
\phi(x)=\psi(x)=0,&x=g(t)\, \textrm{or}\, x=h(t).
\end{array} \right.
\label{B1f}
\end{eqnarray}
\end{lem}

It follows from Theorems 2.2 and \ref{basic2} that
\begin{thm} $R_0^F(t)$ is strictly monotone increasing function of $t$, which means that if $t_1<t_2$, then $R_0^F(t_1)<R_0^F(t_2)$.
Moreover, $\lim_{t\to \infty}\, R_0^F(t)\geq \sqrt{\frac{\frac{r_\infty\gamma\infty}{{\mu_2}_\infty+\gamma_\infty}}{{\mu_1}_\infty+\frac{\nu^2}{4D^2}}}$ if $(H)$ holds and $h(t)-g(t)\to \infty$ as $t\to \infty$.
\end{thm}

\begin{rmk} We have assumed that $(H)$ holds and $\nu<2D\sqrt{\frac{r_\infty\gamma_\infty-{\mu_1}_\infty({\mu_2}_\infty+\gamma_\infty)}{{\mu_2}_\infty+\gamma_\infty}}$ in this paper. By Theorem 3.4, we have
 $R_0^F(t_0)>1$ for some $t_0>0$ provided that $h(t)-g(t)\to \infty$ as $t\to \infty$.\label{rem1}
\end{rmk}

\section{Mosquitos vanishing}
It follows from Theorem 2.2 that $x=g(t)$ is monotonic decreasing and $x=h(t)$ is monotonic increasing, therefore
there exist  $h_\infty, -g_\infty\in (0, +\infty]$ such that  $\lim_{t\to +\infty} \, g(t)$ $=g_\infty$
and $\lim_{t\to +\infty} \ h(t)=h_\infty$. The following lemma shows that  both $g_\infty$ and $h_\infty$ are finite or infinite simultaneously, that is, if $h_\infty<\infty$, then $-g_\infty<\infty$, vice versa.

\begin{lem}
If $h_\infty<\infty$ or $g_\infty>-\infty$, then both $h_\infty$ and
$g_\infty$ are finite and
$$
R_0^{DA}((g_\infty, h_\infty),D, \nu)\leq 1\ \textrm{and}\
\lim_{t\to\infty} (\|M(t, \cdot)\|_{C([g(t),\,h(t)])}+\|A(t, \cdot)\|_{C([g(t),\,h(t)])})=0.$$
\label{spr}
\end{lem}
\bpf Without loss of generality, we assume that $h_\infty<\infty$, and prove that $R_0^{DA}\leq 1$, which implies that $g_\infty>-\infty$ by Remark \ref{rem1}.

{\bf Step 1.}
We first prove that $\lim_{t\rightarrow\infty}h'(t)=0$. In fact, the transformation $y=\frac{x}{h(t)}h_0$, $w(t,y)=M(t,\frac{h(t)}{h_0}y)=M(t,x)$ turns problem \eqref{a3} for $M(t,x)$ in $[0,+\infty)\times[0,h(t)]$ into a new problem for $w(t,y)$ in $[0, +\infty)\times[0, h_0]$. Let $\chi$ be the function in $C^3([0,h_0])$ satisfying
\begin{eqnarray*}
\left\{
\begin{array}{lll}
\chi(y)=1, &
\frac{h_0}{2}\leq y \leq h_0,  \\
\chi(y)=0,\; &0\leq y \leq \frac{h_0}{8}.
\end{array} \right.
\end{eqnarray*}
Letting $z(t,y)=w(t,y)\times \chi(y)$, where $(t,y)\in [0,+\infty)\times[0,h_0]$, it follows from the $L^p$ theory of parabolic equations and the Sobolev imbedding theory that there exists $M_1>0$ such that
$$\|z\|_{C^{(1+\alpha)/2,1+\alpha
}([0, +\infty)\times[0, h_0])}\leq M_1,$$
therefore, there exists $M_2>0$ such that $\|M\|_{C^{(1+\alpha)/2,1+\alpha
}([0, +\infty)\times[0, h(t)]) }\leq M_2$. Lemma 2.3 together with the free boundary condition yields that there exists $M_3>0$ such that
$$\|h\|_{C^{1+\alpha
}([0, +\infty)) }\leq M_3,$$
which together with the assumption $h_\infty<\infty$ implies that $\lim_{t\rightarrow\infty}h'(t)=0$.

Next we
assume that $R_0^{DA}((g_\infty, h_\infty), D, \nu)>1$ by contradiction. Similarly as Lemma 3.1 in \cite{MZ}, we know that there is $\varepsilon_0>0$ such that $h'(t)>\varepsilon_0$. This contradicts the fact $\lim_{t\rightarrow\infty}h'(t)=0$.

{\bf Step 2.} $\lim_{t\to
+\infty} \ ||M(t,\cdot)||_{C([g(t),\, h(t)])}=\lim_{t\to
+\infty} \ ||A(t,\cdot)||_{C([g(t),\, h(t)])}=0$.

Let $(\overline M(t,x), \overline A(t,x))$ denote the unique solution of the problem
{\small\begin{eqnarray}
\left\{
\begin{array}{lll}
\overline M_t-D\overline M_{xx}+\nu \overline M_x=\gamma(x)\overline A(1-\frac{\overline M}{K_1})-\mu_1(x)\overline M,
&t>0,\,g_\infty<x<h_\infty , \\
\overline A_t=r(x)(1-\frac{A}{K_2})-(\mu_2(x)+\gamma (x))\overline A,&t>0,\,g_\infty<x<h_\infty , \\
\overline M(0,g_\infty)=\overline A(0,g_\infty)=\overline M(0,h_\infty)=\overline A(0,h_\infty)=0, & t>0,\\
(\overline M(0,x), \overline A(0,x))=(\tilde M_0(x), \tilde A_0(x)), &g_\infty\leq x\leq h_\infty,
\end{array} \right.
\label{m23}
\end{eqnarray}}
with
\begin{eqnarray*}
(\tilde M_0(x), \tilde A_0(x))= \left\{
\begin{array}{lll}
(M_0(x), A_0(x)),&g_0\leq x\leq h_0, \\
(0,0), & \mbox{ otherwise}.
\end{array} \right.
\end{eqnarray*}
It follows from the comparison principle that $(0, 0)\leq (M, A)(t,x)\leq (\overline M, \overline A)(t, x)$
for $t>0$ and $x\in [g(t), h(t)]$.

Using the fact $R_0^{DA}((g_\infty, h_\infty),D, \nu)\leq 1$ in step 1, we find that
$(0, 0)$ is the unique nonnegative steady-state solution of problem (\ref{m23}). Choosing the lower solution as $(0, 0)$ and upper solution as $(K_1, K_2)$, it is easy to see, by the method of upper and lower solutions and
its associated monotone iterations, that the time-dependent solution converges to the unique nonnegative steady-state solution.
Therefore, $(\overline M(x,t), \overline A(x,t))\to (0,0)$ uniformly
for $x\in [g_\infty, h_\infty]$ as $t\to \infty$ and then $\lim_{t\to
+\infty} \ ||M(t,\cdot)||_{C([g(t),\, h(t)])}=\lim_{t\to
+\infty} \ ||A(t,\cdot)||_{C([g(t),\, h(t)])}=0$.
 \epf
\bigskip

 Next, for that
 the invasive regime of the Aedes aegypti mosquitoes depends on whether  $h_\infty-g_\infty=\infty$ and $\lim_{t\to
+\infty}\ (||M(t, \cdot)||_{C(g(t), h(t)])}+||A(t, \cdot)||_{C([g(t),h(t)])})=0$, we have the following definitions:

\begin{defi}
We say that vanishing occurs or mosquitoes vanish eventually if
$$h_\infty-g_\infty <\infty\, \textrm{ and}\,
 \lim_{t\to +\infty} \, (||M(t, \cdot)||_{C([g(t),h(t)])}+||A(t, \cdot)||_{C([g(t), h(t)])})=0,$$
  and  spreading occurs or mosquitoes spread successfully if $$h_\infty-g_\infty =\infty\, \textrm{and}\,
\limsup_{t\to +\infty}\, (||M(t, \cdot)||_{C([g(t),h(t)])}+||A(t, \cdot)||_{C([g(t),h(t)])})>0.$$
\end{defi}

 The next result shows that if $h_\infty-g_\infty<\infty$, then vanishing occurs.
\begin{lem}  If $h_\infty-g_\infty<\infty$, then $\lim_{t\to
+\infty} \, (||M(t,\cdot)||_{C([g(t),h(t)])}+||A(t,\cdot)||_{C([g(t),h(t)])})=0$.
\end{lem}
\bpf
Assume that
$$\limsup_{t\to +\infty} \ ||M(t,\cdot)||_{C([g(t), h(t)])}=\delta>0$$
 by contradiction. Then there exists a sequence $( t_k,x_k )$
in $(0, \infty)\times(g(t), h(t)) $
such that $M( t_k,x_k)\geq \delta /2$ for all $k \in \mathbb{N}$, and $t_k\to \infty$ as $k\to \infty$.

Now we claim that
\begin{eqnarray}
\|M\|_{C^{(1+\alpha)/2,1+\alpha}}\leq C,\ t\geq 1 ,\label{est-2}\\
||h'||_{C^{\alpha/2}([1, +\infty))},\, ||g'||_{C^{\alpha/2}([1, +\infty))}\leq  C
\label{est-1}
\end{eqnarray}
for any $\alpha \in(0,1)$ and some positive constant $C$. In fact, straighten the double free boundary fronts by the transformation
$$y=\frac{2h_0x}{h(t)-g(t)}-\frac{h_0(h(t)+g(t))}{h(t)-g(t)},$$
let $w(t,y)=M(t,x),$
then the free boundary problem (\ref{a3}) is transformed into the initial boundary problem (\ref{Hb}) in $(-h_0,h_0)$.
Using the fact that $-g(t)$ and $h(t)$ are increasing and bounded, it follows from standard $L^p$ theory and  the Sobolev imbedding
theorem (\cite{LSU, Lie}) that for $0<\alpha <1$,
there exists a constant $C_1$
depending on $\alpha, h_0, \|M_{0}\|_{C^{2}[-h_0, h_0]}$, $\|A_{0}\|_{C^{2}[-h_0, h_0]}$, $g_\infty, h_\infty$ such that
\begin{eqnarray}\|w\|_{C^{(1+\alpha)/2,1+\alpha
}([\tau, \tau+1])\times[-h_0, h_0] }\leq C_1\label{Bg1}
\end{eqnarray}
for any $\tau\geq 1$. Note that $C_1$ is independent of $\tau$, by using the free boundary conditions in (\ref{a3}), it is easy to see that (\ref{est-2}), (\ref{est-1}) hold.
 Using (\ref{est-1}) and the assumption that $h_\infty-g_\infty<\infty$ yields
$$h'(t)\to 0\ \textrm{and}\ g'(t)\to 0\ \textrm{as}\, t\to +\infty.$$
It follows from the free boundary condition that $\frac {\partial M}{\partial x}(t_k,h(t_k))\to 0$ as $t_k\to \infty$.

On the other hand, since  $-\infty<g_\infty<g(t)<x_k<h(t)<h_\infty<\infty$, there exists a subsequence $\{{x_k}_n\}$ which converges
to $x_0\in [g_\infty, h_\infty]$ as $n\rightarrow\infty$. For convenience, we denote $\{{x_k}_n\}$ as $\{x_k\}$,  it follows that $x_k\to x_0 \in [g_\infty, h_\infty]$ as $k\to \infty$. Thanks to the uniform bound in (\ref{est-2}), we can obtain that $x_0\in (g_\infty, h_\infty)$.

Define $W_k(t,x)=M(t_k+t,x)$ and $Z_k(t,x)=A(t_k+t,x)$  for
$x\in (g(t_k+t), h(t_k+t)), t\in (-t_k, \infty)$.
According to the parabolic regularity, $\{(W_k, Z_k)\}$ has a subsequence $\{(W_{k_i}, Z_{k_i})\}$ which converges to $(\tilde W, \tilde Z)$ as $i\to \infty$, and $(\tilde W, \tilde Z)$ satisfies
\begin{eqnarray*} \left\{
\begin{array}{lll}
\tilde W_t-D \tilde W_{xx}=-\nu\tilde W_x+\gamma(x)\tilde Z(1-\frac{\tilde W}{K_1})-\mu_1(x)\tilde W,& \\
\tilde Z_t=r(x)(1-\frac{\tilde Z}{K_2})\tilde W-(\mu_2(x)+\gamma(x))\tilde Z&
\end{array} \right.
\end{eqnarray*}
for $t\in (-\infty, \infty)$, $g_\infty<x<h_\infty$. Since $\tilde W(t,x_0)\geq \delta/2$, we have $\tilde W>0$ in $(-\infty, \infty)\times(g_\infty, h_\infty)$ by the strong comparison principle.
 Using the Hopf lemma at the point $(0,h_\infty)$ yields
 $\tilde W_x(0,h_\infty)\leq -\sigma^*$ for some $\sigma^*>0$.

Furthermore, the fact $\|M\|_{C^{(1+\alpha)/2,1+\alpha}}\leq C$ implies that
$\frac {\partial M}{\partial x}(t_k+0,h(t_k))=(W_k)_x(0,h(t_k))\to \tilde W_x(0,h_\infty)$ as $k\to \infty$, and then $\tilde W_x(0,h_\infty)=0$, which is a contradiction to $\tilde W_x(0,h_\infty)\leq -\sigma^*<0$,
 Thus $\lim_{t\to +\infty} \ ||M(t,\cdot)||_{C([g(t),h(t)])}=0$.

 Note that $A(t,x)$ satisfies
 $$\frac{\partial A(t,x)}{\partial t}=r(x)(1-\frac{A}{K_2})M-(\mu_2(x)+\gamma(x))A, \ t>0,\, g(t)<x<h(t),$$
 and $r(x)(1-\frac{A}{K_2})M\to 0$ uniformly for $x\in [g(t), h(t)]$ as $t\to \infty$, we then have $\lim_{t\to
+\infty} \ ||A(t, \cdot)||_{C([g(t),h(t)])}=0$.
\epf

\bigskip
Now we give sufficient conditions so that the mosquitoes are vanishing.
\begin{lem} If $R_0^F(\infty)\leq 1$, then $h_\infty-g_\infty<\infty$ and
$\lim_{t\to +\infty} \ (||M(t, \cdot)||_{C([g(t),h(t)])}+||A(t, \cdot)||_{C([g(t),h(t)])})=0$.
\label{vanish}
\end{lem}

In this paper we assume that the far site is high-risk and consider small advection, so if $h_\infty-g_\infty=\infty$, then $R_0^F(\infty)> 1$. Therefore $R_0^F(\infty)\leq1$ means that vanishing happens. The next result shows that if $R_0^F(0)<1$, vanishing occurs for small initial values.
\begin{thm} If $R_0^F(0)<1$, then $h_\infty-g_\infty<\infty$ and
$\lim_{t\to +\infty} \ (||M(t, \cdot)||_{C([g(t),h(t)])}+||A(t, \cdot)||_{C([g(t),h(t)])})=0$  provided that $||M_0(x)||_{C([-h_0, h_0])}$, $||A_0(x)||_{C([-h_0, h_0])}$ are sufficiently small.
\end{thm}
\bpf We construct a suitable upper solution to problem \eqref{a3}.
Since $R_0^F(0)<1$, it follows from Lemma 3.3 that there is a $\lambda_0>0$ and $0<\phi(x), \psi(x)\leq 1$ in $(-h_0, h_0)$ such that
\begin{eqnarray}
\left\{
\begin{array}{lll}
-D \phi_{xx}+\nu\phi_x= \gamma(x)\psi-\mu_1(x)\phi+\lambda_0 \phi,\; &
-h_0<x<h_0,  \\
0= r(x)\phi-(\mu_2(x)+\gamma(x))\psi+\lambda_0 \psi,\; &
-h_0<x<h_0,  \\
\phi(x)=\psi(x)=0, &x=\pm h_0.
\end{array} \right.
\label{B1f1}
\end{eqnarray}

Recalling that $\phi(h_0), \psi(h_0)<0$ and $\phi(-h_0), \psi(-h_0)>0$ yields that there exist some positive constants $C_1$ and $C_2$ such that
  $$x\phi'\leq C_1\phi,\  x\psi'\leq C_2\psi\ \textrm{for}\ -h_0<x<h_0.$$
Also, there exists $L>0$ such that
\begin{eqnarray}
\frac{1}{L}\leq\frac{\phi(x)}{\psi(x)}\leq L.
\label{B10f1}
\end{eqnarray}
In fact, for the right inequality in \eqref{B10f1}, it is easy to see that $\phi'(h_0),\psi'(h_0)<0$, thus there exists $\delta_1>0$ such that for $ x\in [h_0-\delta_1, h_0]$,
$$\phi'(x)<\frac{\phi'(h_0)}{2}<0,\ \psi'(x)<\frac{\psi'(h_0)}{2}<0.$$
Setting $L_1=\max_{[h_0-\delta_1,h_0]}\{\frac{\phi'(x)}{\psi'(x)}\}$,
we have $$\frac{\phi'(x)}{\psi'(x)}\leq L_1,$$
as a result, we have $\phi'(x)\geq L_1\psi'(x)$ in $[h_0-\delta_1,h_0]$,
it follows from the mean-value theorem that $$\frac{\phi(x)}{\psi(x)}\leq L_1 \ \mbox{in}\ [h_0-\delta_1,h_0].$$
Similarly, setting $L_2=\max_{[h_0-\delta_1,h_0]}\{\frac{\psi'(x)}{\phi'(x)}\}$, we obtain
 $$\frac{\psi(x)}{\phi(x)}\leq L_2 \ \mbox{in}\ [h_0-\delta_1,h_0].$$
 Let $L_3=\max\{ L_1, L_2\}$, it follows that
$$\frac{1}{L_3}\leq\frac{\phi(x)}{\psi(x)}\leq L_3,\ x\in [h_0-\delta_1, h_0].$$
Similarly, for the left inequality in \eqref{B10f1}, there exist $\delta_2>0$ and $L_4>0$ such that for $x\in [-h_0, -h_0+\delta_2]$,
$$\frac{1}{L_4}\leq\frac{\phi(x)}{\psi(x)}\leq L_4.$$
For the remain part $x\in(-h_0+\delta_2,h_0-\delta_1),$ since $\phi(x),\psi(x)$ are both positive, there exists $L_5>0$, such that $$\frac{1}{L_5}\leq\frac{\phi(x)}{\psi(x)}\leq L_5.$$
Therefore \eqref{B10f1} holds for $L:=\max\{L_3, L_4, L_5\}$.

Similarly as in \cite{DL}, we set
$$\sigma (t)=h_0(1+\delta-\frac \delta 2 e^{-\delta t}), \  t\geq 0,$$
and
$$\overline M=\varepsilon e^{-\delta t}\phi(xh_0/\sigma (t))e^{\frac{\nu}{2D}(1-\frac{h_0}{\sigma(t)})x}, \ -\sigma(t)\leq
x\leq \sigma(t),\ t\geq 0.$$
$$\overline A=\varepsilon e^{-\delta t}\psi(xh_0/\sigma (t))e^{\frac{\nu}{2D}(1-\frac{h_0}{\sigma(t)})x}, \ -\sigma(t)\leq
x\leq \sigma(t),\ t\geq 0.$$

Since $\lambda_0>0$, it follows from  \eqref{B10f1} and the continuity of the function $r(x)$, $\mu_1(x)$, $\mu_2(x)$ and $\gamma(x)$ in $[-2h_0, 2h_0]$ that there exists a small $\delta >0$ such that
\begin{eqnarray*}
-\delta-\frac{\nu h_0}{4D}(1+\delta)\frac{{h_0}^2}{\sigma^2(t)}\delta^2-\frac{{h_0}^2}{\sigma^2(t)}\frac{\delta^2}{2}C_1+\frac{\nu^2}{4D}(1-\frac{{h_0}^2}{\sigma^2(t)})
+\frac{{h_0}^2}{\sigma^2(t)}\lambda_0\\
-L|\frac{{h_0}^2}{\sigma^2(t)}\gamma(y)-\gamma(x)|+(\mu_1(x)-\frac{{h_0}^2}{\sigma^2(t)}\mu_1(y))\geq0,
\end{eqnarray*}
and
\begin{eqnarray*}
-\delta-\frac{\nu h_0}{4D}(1+\delta)\frac{{h_0}^2}{\sigma^2(t)}\delta^2-\frac{{h_0}^2}{\sigma^2(t)}\frac{\delta^2}{2}C_2+\lambda_0-L|r(y)-r(x)|\\
+(\mu_2(x)-\mu_2(y))+(\gamma(x)-\gamma(y))\geq0,
\end{eqnarray*}
where $y=\frac{xh_0}{\sigma(t)}.$

Direct computations yield
\begin{eqnarray*}
& &\overline M_t-D\overline M_{xx}+\nu\overline M_x-\gamma(x)\overline A(1-\frac{\overline M}{K_1})+\mu_1(x)\overline M \\
& &\geq\overline M_t-D\overline M_{xx}+\nu\overline M_x-\gamma(x)\overline A+\mu_1(x)\overline M \\
& &=-\delta\overline M-\frac{\nu x}{2D}\frac{{h_0}^2}{\sigma^2(t)}\frac{\delta^2}{2}e^{-\delta t}\overline M
-\frac{x{h_0}^2}{\sigma^2(t)}\frac{\delta^2}{2}e^{-\delta t}\overline M \phi^{-1}\phi'
+\frac{\nu^2}{4D}(1-\frac{{h_0}^2}{\sigma^2(t)})\overline M\\
& &+\lambda_0\frac{{h_0}^2}{\sigma^2(t)}\overline M+\varepsilon e^{-\delta t}e^{\frac{\nu}{2D}(1-\frac{h_0}{\sigma(t)})x}\psi(\frac{{h_0}^2}{\sigma^2(t)}\gamma(y)-\gamma(x))\\
& &+ \varepsilon e^{-\delta t}e^{\frac{\nu}{2D}(1-\frac{h_0}{\sigma(t)})x}\phi(\mu_1(x)-\frac{{h_0}^2}{\sigma^2(t)}\mu_1(y))    \\
& &\geq \overline M(-\delta-\frac{\nu h_0}{4D}(1+\delta)\frac{{h_0}^2}{\sigma^2(t)}\delta^2-\frac{{h_0}^2}{\sigma^2(t)}\frac{\delta^2}{2}C_1+\frac{\nu^2}{4D}(1-\frac{{h_0}^2}{\sigma^2(t)})\\
& &+\frac{{h_0}^2}{\sigma^2(t)}\lambda_0
-L|\frac{{h_0}^2}{\sigma^2(t)}\gamma(y)-\gamma(x)|+(\mu_1(x)-\frac{{h_0}^2}{\sigma^2(t)}\mu_1(y))) \\
& &\geq 0,
\end{eqnarray*}
\begin{eqnarray*}
& &\overline A_t-r(x)(1-\frac{\overline A}{K_2})\overline M+(\mu_2(x)+\gamma(x))\overline A\\
& &\geq\overline A_t-r(x)\overline M+(\mu_2(x)+\gamma(x))\overline A\\
& &=-\delta\overline A+\frac{\nu x}{2D}\frac{{h_0}^2}{\sigma^2(t)}\frac{\delta^2}{2}e^{-\delta t}\overline A
-\frac{x{h_0}^2}{\sigma^2(t)}\frac{\delta^2}{2}e^{-\delta t}\overline A \psi^{-1}\psi'+(\lambda_0\overline A+r(y)\overline M\\
& &-(\mu_2(y)+\gamma(y))\overline A)-r(x)\overline M+(\mu_2(x)+\gamma(x))\overline A\\
& &\geq\overline A(-\delta-\frac{\nu h_0}{4D}(1+\delta)\frac{{h_0}^2}{\sigma^2(t)}\delta^2-\frac{{h_0}^2}{\sigma^2(t)}\frac{\delta^2}{2}C_2+\lambda_0-L|r(y)-r(x)|\\
& &+(\mu_2(x)-\mu_2(y))+(\gamma(x)-\gamma(y)))\\
& &\geq0,
\end{eqnarray*}
for all $-\sigma (t)<x<\sigma (t)$ and $t>0$.

On the other hand,
we can choose $\varepsilon=\frac {\delta^2h_0} {2\mu e^{\frac{\nu}{2D}h_0\delta}}\min\{\frac{-1}{\phi'(h_0)},\frac{1}{\phi'(-h_0)}\}$ such that
\begin{eqnarray*}
\left\{
\begin{array}{lll}
\overline M_t\geq D\overline M_{xx}+\nu\overline M_x-\gamma(x)\overline A(1-\frac{\overline M}{K_1})+\mu_1(x)\overline M,\; &  -\sigma(t)<x<\sigma(t), \, t>0,\\
\overline A_t\geq r(x)(1-\frac{\overline A}{K_2})\overline M-(\mu_2(x)+\gamma(x))\overline A,\; &  -\sigma(t)<x<\sigma(t), \, t>0,\\
\overline M(t,x)=\overline A(t,x)=0,&x=\pm \sigma(t)\, \, t>0,\\
-\sigma (0)<-h_0,\; -\sigma'(t)\leq -\mu \frac{\partial \overline M}{\partial x}(t, -\sigma(t)), & t>0, \\
\sigma(0)> h_0, \; \sigma'(t)\geq -\mu \frac{\partial \overline M}{\partial x}(t, \sigma(t)), & t>0.
\end{array} \right.
\end{eqnarray*}
If  $||M_{0}||_{L^\infty(-h_0,h_0)}\leq \varepsilon \phi(\frac {h_0}{1+\delta/2})e^{\frac{\nu}{2D}\frac{\delta}{2+\delta}(-h_0)}$
and $||A_{0}||_{L^\infty(-h_0,h_0)}\leq \varepsilon \psi(\frac {h_0}{1+\delta/2})e^{\frac{\nu}{2D}\frac{\delta}{2+\delta}(-h_0)}$, then
$M_{0}(x)\leq \varepsilon \phi(\frac {x}{1+\delta/2})e^{\frac{\nu}{2D}\frac{\delta}{2+\delta}x}=\overline M(0,x)$
and $A_{0}(x)\leq \varepsilon \psi(\frac {x}{1+\delta/2})e^{\frac{\nu}{2D}\frac{\delta}{2+\delta}x}=\overline A(0,x)$
 for $x\in [-h_0, h_0]$, thus $(\overline M(t,x), \overline A(t,x), -\sigma(t), \sigma(t))$ is an upper solution of problem \eqref{a3}.
Applying Lemma 2.3 gives that $g(t)\geq -\sigma(t)$ and $h(t)\leq\sigma(t)$ for $t>0$. It
follows that $h_\infty-g_\infty\leq \lim_{t\to\infty}
2\sigma(t)=2h_0(1+\delta)<\infty$, and $\lim_{t\to +\infty} \ (||M(t, \cdot)||_{C([g(t),h(t)])}+||A(t, \cdot)||_{C([g(t),h(t)])})=0$ by Lemma 4.2.
 \epf

\bigskip
 From the proof above, we have the following result, see Lemma 3.8 in \cite{DL} for details.
\begin{thm} Suppose $R_0^F(0)(:=R_0^{DA}((-h_0, h_0),D, \nu))<1$.  Then $h_\infty-g_\infty<\infty$ and
$$\lim_{t\to +\infty}  (||M(t, \cdot)||_{C([g(t), h(t)])}+||A(t, \cdot)||_{C([g(t), h(t)])})=0$$
 if $\mu$ is sufficiently small.
\end{thm}

\section{Mosquitos spreading}

In this section, we are going to give the sufficient conditions for the mosquitoes to spread. We first
prove that if $R_0^F(0)\geq 1$, the mosquitoes are spreading.
\begin{thm} If $R_0^F(0)\geq 1$, then $h_\infty-g_\infty=\infty$ and $\liminf_{t\to
+\infty} ||M(t, \cdot)||_{C([0, h(t)])}>0$, that is, spreading occurs.
\label{spread}
\end{thm}
\bpf We first consider the case that $R_0^F(0):=R_0^D((-h_0, h_0))>1$. In this case, the linear eigenvalue problem
\begin{eqnarray}
\left\{
\begin{array}{lll}
-d \phi_{xx}+\nu\phi_x= -\gamma(x)\psi-\mu_1(x\phi)+\lambda_0 \phi,\; &
-h_0<x<h_0,  \\
0=r(x)\phi-(\mu_2(x)+\gamma(x))\psi+\lambda_0\psi,\; &
-h_0<x<h_0,  \\
\phi(x)=\psi(x)=0, &x=\pm h_0 .
\end{array} \right.
\label{B2f}
\end{eqnarray}
admits a positive solution $(\phi(x), \psi(x))$ with $||\phi||_{L^\infty}+||\psi||_{L^\infty}=1$, where $\lambda_0$ is the principal eigenvalue. It follows from Lemma 3.3
that $\lambda_0<0$.

We  construct a suitable lower solution to
\eqref{a3} by define
$$\underline {M}(t,x)=\delta \phi(x),\quad \underline A(t,x)=\delta \psi(x)$$
for $-h_0\leq x\leq h_0$, $t\geq 0$, where $\delta $ is sufficiently small such that
\begin{eqnarray*}
\lambda_0+||r||_{L^\infty}\frac{\delta}{K}<0,
\end{eqnarray*}
where $K=\min{\{K_1, K_2\}}.$

 Direct computations yield
\begin{eqnarray*}
& &\underline M_t-D\underline M_{xx}+\nu\underline M_x-\gamma(x)\underline A(1-\frac{\underline M}{K_1})+\mu_1(x)\underline M\\
& &=-D\delta\phi_{xx}+\nu\delta\phi_x-\gamma(x)\delta\psi(1-\frac{\delta\phi}{K_1})+\mu_1(x)\delta\phi\\
& &=\delta\phi(\lambda_0+\gamma(x)\frac{\delta\psi}{K_1})\\
& &\leq0,\\
& &\underline A_t-r(x)(1-\frac{\underline A}{K_2})\underline M+(\mu_2(x)+\gamma(x))\underline A\\
& &=-r(x)\delta\phi+r(x)\frac{\delta\phi}{K_2}\delta\phi+(\mu_2(x)+\gamma(x))\delta\psi\\
& &=\delta\psi(\lambda_0+r(x)\frac{\delta\phi}{K_2})\\
& &\leq0.
\end{eqnarray*}
for all $-h_0<x<h_0$ and $t>0$. Then we have
\begin{eqnarray*}
\left\{
\begin{array}{lll}
\underline M_t\leq D\underline M_{xx}+\nu\underline M_x-\gamma(x)\underline A(1-\frac{\underline M}{K_1})+\mu_1(x)\underline M,\; & t>0 ,\, -h_0<x<h_0 , \\
\underline A_t\leq r(x)(1-\frac{\underline A}{K_2})\underline M+(\mu_2(x)+\gamma(x))\underline A,\; &  t>0 ,\, -h_0<x<h_0 ,\\
\underline M(t, x)=\underline A(t, x)=0,&t>0\, \, x=\pm h_0,\\
0=g'(0)\geq -\mu \underline M_x(t, -h_0), & t>0, \\
0=h'(0)\leq -\mu \underline M_x(t, h_0), & t>0,\\
\underline {M}(0, x)\leq M_{0}(x),\ \underline{A}(0, x)\leq A_{0}(x),\; &-h_0\leq x\leq h_0.
\end{array} \right.
\end{eqnarray*}
Hence, applying Remark 2.1 yields that $M(t, x)\geq\underline M(t, x)$ and  $A(t, x)\geq\underline A(t, x)$
in $[-h_0, h_0]\times [0,\infty)$. It follows that $\liminf_{t\to
+\infty} \ ||M(t, \cdot)||_{C([g(t), h(t)])}\geq \delta \phi(0)>0$, therefore $h_\infty-g_\infty=+\infty$ by Lemma 3.2.

If $R_0^F(0)=1$,  then for any positive time $t_0$, we have $g(t_0)<-h_0$ and $h(t_0)>h_0$; therefore, $R_0^F(t_0)>
R_0^F(0)=1$ by the monotonicity in Lemma 3.4. We then have $h_\infty-g_\infty=+\infty$ as above by replacing the initial time $0$ with the positive time $t_0$.
 \epf

\begin{rmk} It follows from the above proof that spreading occurs, if and only if there exists $t_0\geq 0$ such that $R_0^F(t_0)\geq 1$.
\end{rmk}

Theorems 4.5 and 4.6 show that if $R_0^F(0)<1$, vanishing occurs for small initial scale of mosquitoes or small expanding capability $\mu$, and Theorem 3.5 implies that
if $R^F_0(\infty)\leq 1$, vanishing always occurs for any initial values. The next result
shows that spreading occurs for large expanding capability $\mu$, see similar results and the proofs in \cite{DG,DL}.
 \begin{thm} Suppose that $R_0^F(0)<1$. Then $h_\infty-g_\infty=\infty$ if $\mu$ is sufficiently large.
\end{thm}

\begin{thm} (Sharp threshold) Fixed $h_0$, $M_0$ and $A_0$. There exists $\mu^*\in [0, \infty)$
 such that spreading occurs when $\mu> \mu^*$, and vanishing occurs when $0<\mu\leq \mu^*$.
\end{thm}
\bpf
If $R_0^F(0)\geq 1$, we have $\mu^*=0$, since in this case spreading always happens for $\mu>0$ from Theorem 5.1.

For the remaining case $R_0^F(0)<1$. We define
$$\mu^*:=\sup \{\sigma_0: h_\infty (\mu)-g_\infty(\mu)<\infty \ \textrm{for}\ \mu\in (0,\sigma_0]\}.$$
Theorem 4.5 implies that vanishing happens for all small $\mu>0$, therefore, $\mu^*\in
(0, \infty]$. On the other hand, by Theorem 5.2, it is easy to see that spreading happens for all big $\mu$. Thus we have $\mu^*\in
(0, \infty)$, and spreading happens when $\mu> \mu^*$, vanishing occurs when $0<\mu< \mu^*$ by Corollary 2.6.

We now claim that vanishing happens when $\mu=\mu^*$. Otherwise $h_\infty-g_\infty=\infty$ for
$\mu=\mu^*$. Since $\lim_{t\to \infty}\, R_0^F(t)\geq \sqrt{\frac {\frac{r_{\infty}\gamma_{\infty}}{{\mu_2}_{\infty}+\gamma_{\infty}}}{\frac{\nu^2}{4D^2}+{\mu_1}_{\infty}}}>1$,
 there exists $T_0>0$ such that $R_0^F(T_0):=R_0^{DA}((g(T_0), h(T_0),D, \nu)>1$. By
the continuous dependence of $(M, A, g, h)$ on its initial values, we can find small
$\epsilon>0$ such that the solution of (\ref{a3}) with $\mu=\mu^*-\epsilon$, denoted by
$(M_\epsilon,A_\epsilon,g_\epsilon, h_\epsilon)$, satisfies $R_0^{DA}((g_\epsilon(T_0), h_\epsilon(T_0)), D,\nu)>1$. This implies that spreading
happens for the solution $(M_\epsilon,A_\epsilon, g_\epsilon, h_\epsilon)$, which
contradicts the definition of $\mu^*$. The proof is complete.
 \epf

Next, we consider the asymptotic behavior of the solution to  \eqref{a3} when the spreading occurs.

\begin{thm} Suppose that $h_\infty=-g_\infty=\infty$, then the solution to the
 free boundary problem \eqref{a3} satisfies $\lim_{t\to +\infty} \ (M(t,x),A(t,x))=(M^*(x), A^*(x))$
uniformly in any bounded subset of $(-\infty, \infty)$, where $(M^*(x), A^*(x))$ is the unique bounded positive solution of the following problem
    \begin{eqnarray}
\label{stationary01}
\left\{
\begin{array}{lll}
-DM_{xx}+\nu M_x=\gamma(x)A(1-\frac{M}{K_1})-\mu_1(x)M,&-\infty<x<\infty, \\
0=r(x)(1-\frac{A}{K_2})M-(\mu_2(x)+\gamma(x))A,&-\infty<x<\infty.\\

\end{array} \right.
\end{eqnarray}
\end{thm}
\bpf
(1) The existence and uniqueness of the stationary solution

It is easy to see that problem
(\ref{stationary01}) is equivalent to
\begin{eqnarray}
-DM_{xx}+\nu M_x=\gamma(x)\frac{r(x)M}{r(x)\frac{M}{K_2}+\mu_2(x)+\gamma(x)}(1-\frac{M}{K_1})-\mu_1(x)M
\label{stationary2}
\end{eqnarray}
for$-\infty<x<\infty.$
Since $h_\infty=-g_\infty=+\infty$, it follows from Remark 3.1 that there exists $t_0>0$ such that $R_0^F(t_0)=R_0^{DA}((g(t_0),h(t_0)),D,\nu)>1$,
therefore, for any $L$ with $L\geq L_0:=\max\{-g(t_0), h(t_0)\}$, we consider the problem
 \begin{eqnarray}
\label{stationary3}
\left\{
\begin{array}{lll}
-D{M_L}_{xx}+\nu {M_L}_x=\frac{\gamma(x)r(x)M_L}{r(x)\frac{M_L}{K_2}+\mu_2(x)+\gamma(x)}(1-\frac{M_L}{K_1})-\mu_1(x)M_L,&-L<x<L, \\
 M_L(\pm L)=0.&
\end{array} \right.
\end{eqnarray}

Setting $\tilde{M_L}=K_1, \hat{M_L}=\delta\phi(x)$, where $(\phi(x),\psi(x))$ is the corresponding eigenfunction to the principal eigenvalue $\lambda_0$ of the following problem
\begin{eqnarray}
\left\{
\begin{array}{lll}
-D\phi_{xx}+\nu\phi_x=\gamma(x)\psi-\mu_1(x)\phi+\lambda_0\phi,\; &
x\in (-L, L),  \\
0=r(x)\phi-(\mu_2(x)+\gamma(x))\psi+\lambda_0\psi,\; &
x\in (-L, L),  \\
\phi(x)=\psi(x)=0, &x=\pm L.
\end{array} \right.
\label{B1k0}
\end{eqnarray}
We can choose $\delta$ sufficiently small such that $\tilde{M_L}$ and $\hat{M_L}$ are upper and lower solutions to  problem (\ref{stationary3}). As a result, there exists $M_L$ that solves problem (\ref{stationary3}).

Moreover, for the first equation in (\ref{stationary3}), taking $M_L=e^{\frac{\nu}{2D}x}u$ gives that

$-Du_{xx}=-\frac{\nu^2}{4D}u+\gamma(x)\frac{r(x)u}{r(x)\frac{e^{\frac{\nu}{2D}x}u}{K_2}+\mu_2(x)+\gamma(x)}(1-\frac{e^{\frac{\nu}{2D}x}u}{K_1})-\mu_1(x)u:=f(u)u$\\
It is easy to see that $f(u)$ is decreasing, therefore the positive solution is unique.

Using the comparison principle yields that as $L$ increases to infinity, $M_L$ increases to a positive solution $M^*$ to problem
(\ref{stationary2}).
 The uniqueness of positive solution to problem
(\ref{stationary2}) follows from the similar technique in \cite{DLM}.

(2) The limit superior of the solution

We recall that the comparison principle gives $(M(t,x),A(t,x))\leq (\overline M(t,x), \overline A(t,x))$
for $(t,x)\in (0,\infty)\times(-\infty,\infty) $, where
$(\overline M(t,x), \overline A(t,x))$ is the solution of the problem
    \begin{eqnarray}
\label{pde1}
\left\{
\begin{array}{ll}
\overline M_{t}=D \overline M_{xx}-\nu\overline M_{x}+\gamma (x)\overline A(1-\frac{\overline M}{K_1}) -\mu_1(x)\overline M, &t>0,\; -\infty<x<\infty,   \\
\overline A_{t}=r(x)(1-\frac{\overline A}{K_2})\overline M -(\mu_2(x)+\gamma(x))\overline A, &t>0,\; -\infty<x<\infty,   \\
\overline M(0,x)=K_1,\ \overline A(0,x)=K_2.&
\end{array} \right.
\end{eqnarray}
It is easy to see that $(\overline M(t,x), \overline A(t,x))\leq(\overline M(0,x),\overline A(0,x))$, therefore we have $(\overline M(t+\delta,x), \overline A(t+\delta,x))\leq(\overline M(t,x),\overline A(t,x))$ by comparing the initial conditions, that is, $(\overline M, \overline A)$ is monotone decreasing with respect to $t$ and $\lim_{t\to\infty}
(\overline M, \overline A)= (M^*(x),A^*(x))$ uniformly in any bounded subset of $(-\infty, \infty)$; therefore we deduce
\begin{equation}\label{123}
\limsup_{t\to +\infty} \ (M(t,x),A(t,x))\leq (M^*(x),A^*(x))
\end{equation}
uniformly in any bounded subset of $(-\infty, \infty)$.

(3) The lower bound of the solution for a large time

From step 1, we can deduce that the principal eigenvalue $\lambda_0$ of
\begin{eqnarray}
\left\{
\begin{array}{lll}
-D\phi_{xx}+\nu\phi_x=\gamma(x)\psi-\mu_1(x)\phi+\lambda_0\phi,\; &
x\in (-L_0, L_0),  \\
0=r(x)\phi-(\mu_2(x)+\gamma(x))\psi+\lambda_0\psi,\; &
x\in (-L_0, L_0),  \\
\phi(x)=\psi(x)=0, &x=\pm L_0
\end{array} \right.
\label{B1k}
\end{eqnarray}
satisfies
$$\lambda_0<0.$$
Since $h_\infty=\infty=-g_\infty$,
 for any $L\geq L_0$, there exists $t_L>0$ such that $g(t)\leq -L$ and $h(t)\geq L$ for $t\geq t_L$.

Letting $\underline M=\delta \phi$ and $\underline A=\delta\psi$,  we can choose $\delta$ sufficiently small such that $(\underline M, \underline A)$
satisfies
 \begin{eqnarray*}
\left\{
\begin{array}{lll}
\underline M_{t}-\underline M_{xx}+\nu\underline M_x\leq\gamma(x)\underline A(1-\frac{M}{K_1})-\mu_1(x)\underline M,\; &  t>t_{L_0},\-L_0<x<L_0,  \\
\underline A{t}\leq r(x)(1-\frac{A}{K_2})\underline M-(\mu_2(x)+\gamma(x))\underline A,\; &  t>t_{L_0},\-L_0<x<L_0,  \\
\underline M(t,x)=\underline A(t,x)=0,\ \quad & t>t_{L_0},\-L_0<x<L_0 ,\\
\underline M(t_{L_0},x)\leq M(t_{L_0},x),\ \underline A(t_{L_0},x)\leq A(t_{L_0},x),
 & -L_0\leq x\leq L_0,
\end{array} \right.
\end{eqnarray*}
 which means that $(\underline M, \underline A)$ is a lower solution of $(M,A)$ in $[t_{L_0}, \infty)\times[-L_0, L_0] $.
We then have $(M, A)\geq (\delta \phi, \delta\psi)$ in $[t_{L_0}, \infty)\times[-L_0, L_0] $, which implies that the solution can not decay to zero.

(4) The limit inferior of the solution

We extend $\phi(x)$ to $\phi_{L_0}(x)$ by defining $\phi_{L_0}(x):=\phi(x)$ for $-L_0\leq x\leq L_0$ and $\phi_{L_0}(x):=0$ for $x<-L_0$ or $x>L_0$, and $\psi(x)$ to $\psi_{L_0}(x)$ by defining $\psi_{L_0}(x):=\psi(x)$ for $-L_0\leq x\leq L_0$ and $\psi_{L_0}(x):=0$ for $x<-L_0$ or $x>L_0$.
Now for $L\geq L_0$, $(M,A)$ satisfies
{\small \begin{eqnarray}
\left\{
\begin{array}{lll}
M_{t}=D M_{xx}-\nu M_{x}+\gamma (x)A(1-\frac{M}{K_1}) -\mu_1(x)M,\; &t>t_L,\; g(t)<x<h(t),   \\
A_{t}=r (x)(1-\frac{A}{K_2})M -(\mu_2(x)+\gamma(x))A,\; &t>t_L,\; g(t)<x<h(t),   \\
M(t,x)=A(t,x)=0,\, &t>t_L,\,  x=g(t)\, \textrm{or}\, x=h(t),\\
M(t_L,x)\geq \delta \phi_{L_0},\ A(t_L,x)\geq\delta \psi_{L_0},
 & -L\leq x\leq L,
\end{array} \right.
\label{fs1}
\end{eqnarray}}
therefore, we  have $(M, A)\geq (w,z)$ in $[t_L, \infty)\times[-L, L] $,
where $(w,z)$ satisfies
\begin{eqnarray}
\left\{
\begin{array}{lll}
w_{t}=D w_{xx}-\nu w_{x}+\gamma (x)z(1-\frac{w}{K_1}) -\mu_1(x)w,\; &t>t_L,\; -L<x<L,   \\
z_{t}=r (x)(1-\frac{z}{K_2})w -(\mu_2(x)+\gamma(x))z,\; &t>t_L,\; -L<x<L,   \\
w(t,x)=z(t,x)=0,\, &t>t_L,\,  x=\pm L,\\
w(t_L,x)=\delta \phi_{L_0},\ z(t_L,x)=\delta \psi_{L_0},
 & -L\leq x\leq L.
\end{array} \right.
\label{fs11}
\end{eqnarray}
System (\ref{fs11}) is quasimonotone increasing; therefore,
it follows from the upper and lower solution method
 and the theory of monotone dynamical systems ( \cite{HS} Corollary 3.6) that
$\lim_{t\to +\infty} \ (w(t,x), z(t,x))= (M_L(x), A_L(x))\geq(\delta \phi_{L_0},\delta \psi_{L_0})$ uniformly in
$[-L, L]$, where $(M_L(x), A_L(x))$ satisfies
 \begin{eqnarray}
\label{stationary}
\left\{
\begin{array}{lll}
&-DM_{xx}+\nu M_x=\gamma(x)A(1-\frac{M}{K_1})-\mu_1(x)M,&-L<x<L, \\
&0=r(x)(1-\frac{A}{K_2})M-(\mu_1(x)+\gamma(x))A,&-L<x<L.\\
&\quad M(\pm L)=0.\\
\end{array} \right.
\end{eqnarray}

Now we claim the monotonicity and show that if $0<L_1<L_2$, then $(M_{L_1}(x),A_{L_1}(x))$ $\leq$ $(M_{L_2}(x),A_{L_2}(x))$
 in $[-L_1, L_1]$. The result is derived by comparing the initial conditions and boundary
conditions in (\ref{fs11}) for $L=L_1$ and $L=L_2$.

Letting $L\to \infty$, by classical elliptic regularity theory and a diagonal procedure, we obtain that $(M_{L}(x), A_{L}(x))$ converges
uniformly on any compact subset of $(-\infty, \infty)$ to $(M_\infty(x), A_\infty(x))$, which is continuous on $(-\infty, \infty)$ and satisfies
\begin{eqnarray} \left\{
\begin{array}{lll}
-D M''_{\infty}+\nu M'_{\infty}=\gamma(x)A_{\infty}(1-\frac{M_{\infty}}{K_1})-\mu_1(x)M_{\infty},\; &  -\infty<x<\infty,  \\
0=r(x)M_{\infty}(1-\frac{A_{\infty}}{K_2})-(\mu_2(x)+\gamma(x))A_{\infty},\; & -\infty<x<\infty,  \\
M_\infty(x)\geq \delta\phi_{L_0},A_\infty(x)\geq \delta \psi_{L_0}, &-\infty<x<\infty.
\label{fs111}
\end{array} \right.
\end{eqnarray}
It follows from step 1 that $M_\infty (x)= M^*(x)$ and $A_\infty (x)= A^*(x)$.

Now for any given $[-X, X]$ with $X\geq L_0$, since $(M_{L}(x), A_{L}(x))\to (M^*(x), A^*(x))$ uniformly in $[-X, X]$, which is the compact subset of $(-\infty, \infty)$, as $L\to \infty$, we deduce that for any $\varepsilon >0$, there exists $L^*>L_0$ such that  $(M_{L^*}(x), A_{L^*}(x))\geq  (M^*(x)-\varepsilon, A^*(x)-\varepsilon)$ in $[-X, X]$. As above, there is $t_{L^*}$ such that $[g(t), h(t)]\supseteq [-L^*, L^*]$ for $t\geq t_{L^*}$.
Therefore,
$$(M(t,x), A(t,x))\geq (w(t,x), z(t,x))\ \textrm{in}\ [t_{L^*}, \infty)\times [-L^*, L^*],$$ and
$$\lim_{t\to +\infty} \ (w(t,x), z(t,x))= (M_{L^*}(x), A_{L^*}(x))\ \textrm{in}\ [-L^*, L^*].$$
Using the fact that $(M_{L^*}(x),A_{L^*}(x))\geq (M^*(x)-\varepsilon, A^*(x)-\varepsilon)$ in $[-X, X]$ gives
 $$\liminf_{t\to +\infty} \ (M(t,x), A(t,x))\geq (M^*(x)-\varepsilon, A^*(x)-\varepsilon)\ \textrm{in}\ [-X, X].$$
 Since $\varepsilon>0$ is arbitrary, we have $\liminf_{t\to +\infty}  M(t,x)\geq M^*(x)$ and $\liminf_{t\to +\infty} A(t,x)\geq A^*(x)$ uniformly in $[-X,X]$, which together with (\ref{123}) imply that $\lim_{t\to +\infty} \ M(t,x)=M^*(x)$
and  $\lim_{t\to +\infty} \ A(t,x)=A^*(x)$ uniformly in any bounded subset of $(-\infty, \infty)$.
\epf

Combining Remarks 3.1 and 5.1, Lemma 4.1 and Theorem 5.4, we immediately obtain the following
spreading-vanishing dichotomy:
\begin{thm} Suppose that $(H)$ holds and $\nu<2D\sqrt{\frac{r_\infty\gamma_\infty-{\mu_1}_\infty({\mu_2}_\infty+\gamma_\infty)}{{\mu_2}_\infty+\gamma_\infty}}$ .
Let $(M(t, x), A(t, x); g(t), h(t))$ be the solution of  free boundary problem \eqref{a3}.
Then, the following alternatives hold:

Either
\begin{itemize}
\item[$(i)$] {\rm Spreading:} $h_\infty-g_\infty =+\infty$ and $\lim_{t\to +\infty} \ (M(t, x), A(t, x))=(M^*(x), A^*(x))$
uniformly in any bounded subset of $(-\infty, \infty)$; \end{itemize}

or
\begin{itemize}
\item[$(ii)$] {\rm Vanishing:} $h_\infty -g_\infty \leq \infty$ with $R_0^{DA}((g_\infty, h_\infty), D,\nu)\leq1$ and\\
 $\lim_{t\to +\infty} \ (||M(t, \cdot)||_{C([g(t),h(t)])}+||A(t,
\cdot)||_{C([g(t),h(t)])})=0$.
\end{itemize}
\end{thm}

\section{Discussion}

In this paper we constructed a reaction-diffusion-advection model with free boundaries describing the spatial dispersal of Aedes aegypti mosquitoes, which are divided into two life stages. We have obtained some analytical results about the invasive dynamics of Aedes aegypti mosquitoes. The introduction of the thresholds $R_0^{DA}$  for the reaction-diffusion-advection problem with Dirichlet boundary condition and $R_0^F(t)$ for the problem with the free boundary is one of our main results. It is shown that if $R_0^F(t_0)\geq 1$ for some $t_0\geq 0$, mosquitoes spreads successfully (Theorem 5.1, Remark 5.1). If $R_0^F(0)<1$, vanishing happens provided that the initial values of the mosquitoes are small (Theorem 4.4) or the expanding capability is small (Theorem 4.5), while spreading happens for the large expanding capability (Theorem 5.2). Another consideration of our work is the environmental heterogeneity. The complexity of the ecosystem leads to the difference of the habitats mosquitoes survive in. Therefore, the spatial-dependent rates considered in our model confirm more to the reality.

During the outbreaks of mosquito-borne diseases, an emergency measure to reduce the population of Aedes aegypti mosquitoes is insecticides spraying. Minimizing the population of Aedes aegypti mosquitoes is one of the most effective method to control mosquito-borne viruses, such as dengue, Zika, etc. Hence understanding the spatial dispersal dynamics of Aedes aegypti mosquitoes is of great importance. In our model, the free boundary indicates that the spreading or vanishing of the mosquitoes depends on the heterogeneity of the habitats, which have something to do with advection besides the factors discussed above. In this paper we assumed small advection and presented the spreading-vanishing dichotomy, big advection, we believe, will cause more complex transmission dynamics and deserves further study.



%



\begin{thebibliography}{99}\setlength{\itemsep}{-1ex}
{\small

\bibitem{ABL}
I. Ahn, S. Baek, Z. G. Lin, The spreading fronts of an infective environment in a man-environment-man epidemic model, {\it Appl. Math. Model.}, {\bf 40} (2016), 7082¨C7101.

\bibitem{AL} L. J. S. Allen, B. M. Bolker, Y. Lou, A. L. Nevai,
Asymptotic profiles of the steady states for an SIS epidemic reaction-diffusion
model, {\it Discrete Contin. Dyn. Syst. Ser. A} {\bf 21}(2008), 1-20.


\bibitem{CS}
L. Caffarelli and S. Salsa, A Geometric Approach to Free Boundary
Problems, Graduate Studies in Mathematics, {\bf 68}, American
Mathematical Society, Providence, RI, (2005).


\bibitem{CC}
R. S. Cantrell and C. Cosner, Spatial Ecology via Reaction-Diffusion Equations, John Wiley
and Sons Ltd., Chichester, UK, 2003.

\bibitem{CM}
V. Capasso, L. Maddalena,
Convergence to equilibrium states for a reaction-diffusion system modeling the spatial
spread of a class of bacterial and viral diseases,
{\it J. Math. Biol.}, {\bf 13} (1981), 173-184.

\bibitem{SMC}
S. Cauchemez, M. Ledrans, C. Poletto, P.Quenel, H. De Valk, V. Colizza, P. Y. Bo$\ddot{e}$lle, Local and regional spread of chikungunya fever in the Americas, Euro surveill.  {\it Biometrika}, {\bf19} (2014), 20854.

\bibitem{DPW}
W. W. Ding, R. Peng, L. Wei, The diffusive logistic model with a free boundary in a heterogeneous time-periodic environment,
{\it  J. Differential Equations}, {\bf 263} (2017), 2736-2779.

\bibitem{DG}
Y. H. Du and Z. M. Guo, Spreading-vanishing dichotomy in the diffusive
logistic model with a free boundary, II, {\it J. Differential Equations},
{\bf 250} (2011), 4336-4366.

\bibitem{DGP}
Y. H. Du, Z. M. Guo and R. Peng, A diffusive logistic model with a free boundary in time-periodic environment,
 {\it J. Funct. Anal.} {\bf 265} (2013), 2089-2142.

\bibitem{DLiu}
 Y. H. Du, L. S. Li, Remarks on the uniqueness problem for the logistic equation on the entire space, {\it Bull. Austral. Math. Soc.}
 {\bf 73} (2006), 129-137.


\bibitem{DL}
Y. H. Du, Z. G. Lin, Spreading-vanishing dichotomy in the diffusive
logistic model with a free boundary, {\it SIAM J. Math. Anal.} {\bf
42} (2010), 377-405.

\bibitem{DL2}
 Y. H. Du and Z. G. Lin, The diffusive competition model with a free boundary: invasion of a superior or inferior competitor,
{\it Discrete Contin. Dyn. Syst. Ser. B}, {\bf 19} (2014), 3105-3132.

\bibitem{DB}
Y. H. Du, B. D. Lou, Spreading and vanishing in nonlinear diffusion problems with free boundaries,  {\it J. Eur. Math. Soc.}, {\bf 17}
(2015), 2673-2724.

\bibitem{DLM}
 Y. H. Du, L. Ma, Logistic type equations on $R^N$ by a squeezing method involving boundary blow-up solutions, {\it J. London Math. Soc.}
 {\bf 64} (2001), 107-124.
 
 \bibitem{AFDM} A. S. Fauci, D. M. Morens, Zika virus in the Americas - yet another arbovirus threat, New Eng. J. Med. {\bf 374} (2016), 601-604.

\bibitem{DGYL}
D. Gao, Y. Lou, D. He, T. C. Porco, Y. Kuang, G. Chowell, S. Ruan, Prevention and control of Zika as a mosquito-borne and sexually transmitted disease: A mathematical modeling analysis, Sci. Rep.  {\bf 6} (2016) 28070; doi: 10.1038/srep 28070.

\bibitem{JKZH} J. Ge, K. I. Kim, Z. G. Lin, and H. P. Zhu, A SIS reaction-diffusion-advection model in a low-risk
and high-risk domain, {\it  J. Differential Equations}, {\bf259} (2015), 5486-5509.

\bibitem{GW}
J. S. Guo, C. H. Wu, On a free boundary problem for a two-species weak competition system,
{\it J. Dynam. Differential Equations}, {\bf 24} (2012), 873-895.



\bibitem{KY} Y. Kaneko and Y. Yamada, A free boundary problem for a reaction-diffusion equation
appearing in ecology, {\it Adv. Math. Sci. Appl.}, {\bf 259} (2015), 5486-5509.

\bibitem{LSU} O. A. Ladyzenskaja, V. A. Solonnikov and N. N.
Ural'ceva, Linear and Quasilinear Equations of Parabolic Type, Amer.
Math. Soc, Providence, RI, 1968.


\bibitem{LLZ} C. X. Lei, Z. G. Lin, Q. Y. Zhang, The spreading front of invasive species in favorable habitat or unfavorable habitat,
 {\it J. Differential Equations}, {\bf 257} (2014), 145-166.



\bibitem{MZ} M. Li, Z. G. Lin, The spreading fronts in a mutualistic model with advection, {\it Discrete Contin. Dyn. Syst. Ser. B}, {\bf 20} (2015), 2089-2105.

\bibitem{LPW}
H. C. Li, R. Peng, F. B. Wang, Varying total population enhances disease persistence: Qualitative analysis on a diffusive SIS epidemic model, {\it J. Differential Equations}, {\bf 262} (2017), 885-913.

\bibitem{Lie} G. M. Lieberman, Second Order Parabolic Differential Equations, World Scientific Publishing Co. Inc., River Edge, NJ, 1996.


 \bibitem{PY} R. Peng, F. Q. Yi, Asymptotic profile of the positive steady state for an SIS epidemic reaction-diffusion model: effects of epidemic risk and population movement, {\it Phys. D}, {\bf 259} (2013), 8-25.


\bibitem{PZ}
R. Peng, X. Q. Zhao, The diffusive logistic model with a free boundary and seasonal succession,
{\it Discrete Contin. Dyn. Syst. A} (5) {\bf 33} (2013), 2007-2031.

\bibitem{HS}
H. L. Smith, Monotone Dynamical Systems, American Math. Soc., Providence, 1995.

\bibitem{LNW}
L. T. Takahashi, N. A. Maidana, W. C. Ferreira Jr, P. Pulino, H M. Yang, Mathematical models for the Aedes aegypti dispersal dynamics: travelling waves by wing and wind, {\it
Bull. Math. Biol}, {\bf 67} (2005), 509-528.

\bibitem{CTSR}
C. R. Tian, S. G. Ruan, A free boundary problem for Aedes aegypti mosquito invasion, {\it Appl. Math. Model.}, {\bf 46} (2017), 203-217. 



\bibitem{WW2} Y. X. Wang, Z. C. Wang, Entire solutions in a time-delayed and diffusive epidemic model,
{\it Appl. Math. Comput.}, {\bf 219} (2013), 5033-5041.

\bibitem{Wa}
M. X. Wang, On some free boundary problems of the prey-predator model, {\it J. Differential Equations}, {\bf 256} (2014), 3365-3394.


\bibitem{WZ}
M. X. Wang and J. F. Zhao, Free boundary problems for a Lotka-Volterra competition system, {\it J. Dynam. Differential Equations} {\bf 26} (2014),
 655-672.

\bibitem{WS}
S. L. Wu, Entire solutions in a bistable reaction-diffusion system modeling man-environment-man epidemics,
{\it Nonlinear Anal. Real World Appl.}, {\bf 13}(2012), 1991-2005.

\bibitem{ZW} X. Q. Zhao, W. Wang, Fisher waves in an epidemic model, {\it Discrete Contin. Dyn. Syst. B}, {\bf 4} (2004), 1117-1128.

\bibitem{ZH} J. Zhou, H. W. Hethcote, Population size dependent incidence in
models for diseases without immunity, {\it J. Math. Bio.} {\bf 32} (1994), 809-834.


\bibitem{ZX}
P. Zhou and D. M. Xiao, The diffusive logistic model with a free boundary in heterogeneous environment,
{\it J. Differential Equations} {\bf 256} (2014), 1927-1954.

\bibitem{WHO} WHO, Dengue and Dengue hemorrhagic fever, World Health Organization, Ginebra, Suiza, 2002.
 }
\end{thebibliography}
\end{document}